\documentstyle[12pt]{article}
\newcommand\CC{\hbox{C\kern -.58em {\raise .54ex \hbox{$\scriptscriptstyle |$}}
  \kern-.55em {\raise .53ex \hbox{$\scriptscriptstyle |$}} }}
\newcommand\qed{\hfill$\sqcap\kern-8.0pt\hbox{$\sqcup$}$}

\newcommand\N{\hbox{I\kern-.2em\hbox{N}}}
\newcommand\RR{\hbox{I\kern-.2em\hbox{R}}}
\newcommand\sRR{{\sl \hbox{I\kern-.2em\hbox{R}}}}
\newcommand\QQ{\hbox{I\kern-.53em\hbox{Q}}}
\newcommand\ZZ{{{\rm Z}\kern-.28em{\rm Z}}}

\newtheorem{teo}{Theorem}[section]
\newtheorem{prop}[teo]{Proposition}
\newtheorem{lem}[teo]{Lemma}

\newtheorem{coro}[teo]{Corollary}
\newtheorem{df}[teo]{Definition}

\def\fnote#1{\footnote}

\def\abstract{\if@twocolumn
\section*{Abstract}
\else \small {\bf Abstract\vspace{-.5em}\vspace{0pt}} \quotation
\fi}

\begin{document}

\title{\LARGE{Fix-finite approximation property in $F$--spaces}}
\author{\small{\small{Abdelkader Stouti}}}
\date{}
\maketitle

\bigskip

 \begin{center}
\begin{abstract}
In this paper,  with the aid of the simplicial approximation property, the Hopf's construction and Dugundji's homotopy extension Theorem, we first show that if $C$ is a nonempty compact convex  subset of an $F$--space $(E, || \quad ||),$ then  for every $\varepsilon > 0$ and every   subset $D$ of $E$ containing $C$ and every continuous map $f: D\rightarrow C$ there exists a continuous map $g: D\rightarrow C$ which is $\varepsilon$--near to
$f$ and has only a finite number of fixed points. Secondly, by using this result and the
simplicial approximation property, we establish that
for any $\varepsilon > 0$ and every  path and simply connected compact subset $D$ of $E$ containing
$C$ and for each continuous $n$--valued multifunction $F: D\rightarrow 2^{C}$ there exists a
continuous $n$--valued multifunction $G: D\rightarrow 2^{C}$ which is $\varepsilon$--near to
$F$ and has only a finite number of fixed points.
\end{abstract}
\end{center}

\bigskip

\medskip


\medskip

\noindent
\baselineskip=20pt

%%%%%%%%%%%%

\section{Introduction}
Let $(X, d)$ be a metric space, $A$ and $D$ two nonempty subsets of $X$ such that $A \subset D.$  We say that the pair $(D, A)$ enjoy the fix--finite approximation property (F.F.A.P) for a family $\cal F$ of  maps $f : D \rightarrow A$ or multifunctions $F : D \rightarrow 2^{A}$ (where  $2^{A}$ is the set of all nonempty subsets of $A$) if for every
$\varepsilon > 0$ there exists $g \in \cal F$  which  is $\varepsilon$--near to $f$  and has only a finite number of fixed points (respectively there exists a multifunction $G : D \rightarrow 2^{A}$ in $\cal F$ which is $\varepsilon$--near to $F$ and has only a finite number of fixed points). During the last century several authors studies the fix--finite approximation property: Hopf [8], Baillon and Rallis [1] and Schirmer [11].  Later on, in [11, 12] we established some results concerning the fix--finite approximation property in the cases of normed vector spaces and metrizable locally convex vector spaces. In these cases  we have used as key result the Schauder mapping (for examples see: [5, 12, 13]).  In the present paper, we cannot use the Schauder mapping because in $F$--spaces we haven't the useful property of partition of the unity. However,  in this paper we can use the simplicial approximation property which was proved  by Dobrowolski in [6] for  compact convex sets in $F$--spaces.

 In the present paper, we consider the more general case of $F$--spaces. Our key result in this paper is the following:  if $C$ is a nonempty compact convex  subset of an $F$--space $(E, || \quad ||),$ then  for every $\varepsilon > 0$ and every   subset $D$ of $E$ containing $C$ and every continuous map $f: D\rightarrow C$ there exists a continuous map $g: D\rightarrow C$ which is $\varepsilon$--near to
$f$ and has only a finite number of fixed points (Theorem 3.1). By using this result and the
simplicial approximation property, we establish: if $C$ is a nonempty compact convex  subset of an $F$--space $(E, || \quad ||),$ then  for every $\varepsilon > 0$ and for every
 and every  path and simply connected compact subset $D$ of $E$ containing
$C$ and for every continuous $n$--valued multifunction $F: D\rightarrow 2^{C}$ there exists a
continuous $n$--valued multifunction $G: D\rightarrow 2^{C}$ which is $\varepsilon$--near to
$F$ and has only a finite number of fixed points (Theorem 4.1).

\section{Preliminaries}
In this section, we shall recall some definitions and well-know results for subsequent use.

 \bigskip

 Let $\varepsilon > 0$ and let $X$ be a topological space and  $(Y, d)$ be a metric space. We say that $f : X \rightarrow Y$ is a map if it is a single valued function.

\bigskip

\noindent Two continuous maps $f$ and $g$ from $X$ to $Y$ are said to be $\varepsilon$-near if $$ d( f(x), g(x)) < \varepsilon, \hbox{ for all }x \in X. $$

\bigskip

\noindent A homotopy $h_{t} : X \rightarrow Y,\left( 0 \leq t \leq 1 \right)$ is said to be an  $\varepsilon$-homotopy if  $$\sup\{ d( h_{t}(x), h_{t'}(x) ): t, t' \in [0,1]\}< \varepsilon, \hbox{ for all } x \in X. $$

\bigskip

\noindent Two continuous maps $f$ and $g$ from $X$ to $Y$ are said to be $\varepsilon$-homotopic if there exists an $\varepsilon$-homotopy $(h_{t})_{ t \in \left[0,1\right]}$ from $X$ to $Y$ such that $h_{0}=f$ and $h_{1}=g.$

\bigskip

 Let $X$  be a Hausdorff topological space and $f : X \rightarrow X$ be a map. A point $x$ of $X$ is said to be a fixed point of $f$ if $f(x)=x$. We denote by $Fix(f)$ the set of all fixed points of $f.$

\bigskip

 Let $Y$ be a metric space. One says that $Y$ is an absolute neighborhood retract (ANR) if for any nonempty closed subset $A$ of an arbitrary metric space $X$ and for any continuous map $f : A \rightarrow Y,$ then there exists an open subset $U$ of $X$ containing $A$ and a continuous map $g : U \rightarrow Y$ which is an extension of $f$ (i.e. $g(x)=f(x),$ for all  $x \in A$).

\bigskip

\noindent In [7], Dugundji established the homotopy extension Theorem for ANR's.

\begin{teo}
Let $X$ be a metrizable space and $Y$ a ANR. For $\varepsilon>0,$
 there exists $\delta>0$ such that for any two $\delta$-near maps
$f, g : X \rightarrow Y$ and $\delta$-homotopy $j_{t}:A \rightarrow Y,$
where $A$ is a closed subspace of $X$ and $j_{0}=f_{|_{A}}$, $j_{1}=g_{|_{A}}$,
 there exists an $\varepsilon$-homotopy $h_{t} : X \rightarrow Y$ such that
 $ h_{0}=f$, $h_{1}=g$ and
$h_{t_{|_{A}}}=j_{t}$ for all $t \in [ 0,1 ].$
\end{teo}

\bigskip

 Let $X$ and $Y$ be two Hausdorff topological spaces. A multifunction
$F : X \rightarrow  Y$ is a map from $X$ into the  set $2^{Y}$ of nonempty
subsets of $Y.$ The range of $F$ is $F(X)=\displaystyle\bigcup_{x \in X }F(x).$

\bigskip

\noindent The multifunction $F : X \rightarrow  Y$ is said to be upper semi-continuous (usc) if for each open subset $V$ of $Y$ with $F(x) \subset V$ there exists an open subset $U$ of $X$ with $x \in U$ and $F(U)\subset V.$

\bigskip

\noindent The multifunction $F : X \rightarrow  Y$ is called lower semi-continuous (lsc) if for every $x\in X$ and open subset $V$ of $Y$ with $F(x) \cap V \not=\emptyset$ there exists an open subset $U$ of $X$ with $x \in U$ and $F(x')\cap V \not=\emptyset$ for all $x' \in U.$

\bigskip

\noindent The multifunction $F : X \rightarrow  Y$ is continuous if it is both upper semi-continuous and lower semi-continuous. The multifunction $F$ is compact if it is continuous and the  closure of its range $\overline{F(X)}$ is a compact subset of $Y.$

\bigskip

\noindent A point $x$ of $X$ is said to be a fixed point of a multifunction $F : X \rightarrow  X$ if $x \in F(x)$. We denote by $Fix(F)$ the set of all fixed points of $F.$

Let $(X, d)$  be  a metric space. We denote by $C(X)$ the set of nonempty compact subsets of $X.$ Let $A$ and $B$ be two elements of $C(X).$  The Hausdorff distance between  $A$ and $B,$ $d_{H}(A, B),$ is defined by setting:
$$d_{H}(A, B)=\max \left\{\rho(A, B), \rho(B, A) \right\}$$ where
$$\rho(A, B)=\sup\left\{ d(x,B) : x\in A \right\},$$
$$\rho(B, A)=\sup\left\{ d(y,A) : y\in B \right\}$$ and
$$d(x,B)=\inf\left\{ d(x, y) : y\in B \right\}.$$

Let $(X, \tau)$ be a Hausdorff topological space and $(Y, d)$ be a metric space.
Let $F$ and $G$ be two compact multifunctions from $X$ to $Y.$ We define the Hausdorff distance between $F$ and $G$ by setting:
$$d_{H}(F, G)=\sup\left\{ d_{H}(F(x), G(x)) : x \in X \right\}.$$
Let $\varepsilon > 0$ and $F$ and $G$ be two compact multifunctions from $X$ to $Y.$ We say that $F$ and $G$ are $\varepsilon$-near if $d_{H}(F, G)< \varepsilon.$

\bigskip

A multifunction $F : X \rightarrow  X$ is said to be an $n$-function if there exist $n$ continuous maps $f_{i} : X \rightarrow X$, where $i=1,...,n$, such that $F(x)=\{ f_{1}(x),...,  f_{n}(x)\}$ and $f_{i}(x) \not=f_{j}(x)$ for all $x \in X$ and $i,j=1,...,n$ with $i \not=j.$

\begin{df}
Let $X$ and $Y$ be two Hausdorff topological spaces. A multifunction
$F : X \rightarrow Y$ is said to be $n$-valued if for all $x \in X,$ the subset $F(x)$ of $Y$ consists of $n$ points.
\end{df} \par

 Let $X$ and $Y$ be two Hausdorff topological spaces and let $F : X \rightarrow Y$ be a $n$-valued continuous multifunction. Then, we can write $F(x)=\left\{y_{1},..., y_{n}\right\}$ for all $x \in X.$ We define a real function $\gamma$ on  $X$ by
$$\gamma(x)=\min\left\{ || y_{i} - y_{j} || : y_{i}, y_{j}\in F(x), i,j=1,...,n, i \not=j \right\},\hbox{  for all } x \in X.$$
Then,  the gap of $F$ is defined by $\gamma(F)=\inf \left\{ \gamma(x) : x \in X \right\}.$
From [13, p.76], the function $\gamma$ is  continuous.  So if $X$ is compact, then $\gamma(F)> 0.$

\bigskip

 In this paper, we shall need the following result due to H. Schirmer [11].
\begin{lem}$[11].$
Let $X$ and $Y$ be two compact Hausdorff topological spaces. If $X$ is path and simply connected and $F : X\rightarrow Y$ is a continuous $n$-valued multifunction, then $F$ is an $n$-function.
\end{lem}

Next we recall the definition of an $F$--norm.
\begin{df} $[10].$
Let $E$ be a real or a complex linear vector space. A non-negative function $|| \quad
|| : E \rightarrow \RR$ is said to be an $F$--norm if its satisfies the following properties:

(1) $\forall x, y \in E,$  $||x + y|| \leq ||x|| + ||y||$ ;

(2) $||x||=0 \Leftrightarrow x=o_{E}$ ;

(3) $\forall x\in E$ and $\forall t \in [-1, 1],$ $||tx|| \leq ||x||$ ;

(4) if $\alpha_{n} \rightarrow 0,$ then for every $x \in E$ $\alpha_{n}x \rightarrow o_{E}$ ;

(5) the metric $d(x, y)=||x - y||$ is complete.

\end{df}

Throughout this paper we assume that $E$ is a topological vector
space that is not necessarily locally convex. A  linear vector space
equipped with an $F$--norm is called  an $F$-space.

\bigskip

In [10] Kalton and al introduced the simplicial approximation
property which is a useful tool for our results in this paper.

\begin{df} $[10].$
A convex subset $C$ of a metric linear space $(E, || \quad ||)$
has the simplicial approximation property if, for every
$\varepsilon > 0,$ there exists a finite-dimensional compact
convex $C_{\varepsilon} \subset C$ such that, if $S$ is any
finite-dimensional simplex in $C,$ then there exists a continuous
map $h : S \rightarrow C_{\varepsilon}$ with $||h(x) - x || <
\varepsilon$ for every $x \in S.$
\end{df}

Recently in [6, Lemma 2.2] and [6, Corollary 2.6] Dobrowolski
proved the following interesting result.

\begin{lem} $[6].$
Every compact convex set in a metric linear space has the
simplicial approximation property.
\end{lem}

An equivalent formulation of the simplicial approximation property
is given in [10, Theorem 9.8].

\begin{lem} $[10].$
If $K$ is an infinite-dimensional compact convex set in an
$F$-space $E = (E,|| \quad ||)$, then the following statements are
equivalent:

(1) $K$ has the simplicial approximation property,

(2) if $\varepsilon > 0,$ there exist a simplex $S_{\varepsilon}$ in $K$ and a
continuous map $h_{\varepsilon} : K \rightarrow S_{\varepsilon}$ such that $||h_{\varepsilon}(x) - x || <
\varepsilon$ for every $x \in K.$
\end{lem}

\section{The key result}
In this section with the aid of the simplicial approximation property (Lemma 2.7), the
Hopf's construction [3, 8] and Dugundji's homotopy extension Theorem [7]
we shall give our key result in this paper.

 \begin{teo}
 Let  $C$ be a nonempty compact convex  subset of an $F$-space
  $(E, || \quad ||),$ $D$ a subset of $E$ containing $C,$
 $f : D \rightarrow C$ be a continuous map   and $\varepsilon > 0$
be given. Then  there exists a continuous map $g :
D \rightarrow C$ which is $\varepsilon$--near to $f$ and has only
a finite number of fixed points.
 \end{teo}

\bigskip

\noindent {\bf Proof.} Let  $C$ be a nonempty compact convex  subset of an $F$-space
 $(E, || \quad ||),$ $D$ a subset of $E$ containing $C,$
 $f : D \rightarrow C$ be a continuous map and $\varepsilon > 0$ be given.
 From Lemma 2.7, there exist a simplex $S_{\varepsilon}$ in $C$
and  a continuous map $h_{\varepsilon} : C \rightarrow
S_{\varepsilon}$ such that $||h_{\varepsilon}(x) - x || <
\frac{\varepsilon}{2}$ for every $x \in C.$  Put $f_{\varepsilon}=
h_{\varepsilon} \circ f: D \rightarrow  S_{\varepsilon}.$ Hence, $f_{\varepsilon}$ is
$\frac{\varepsilon}{2}$--near to  $f.$ As
$S_{\varepsilon}$ is a simplex, then it is a compact ANR [2].
From Theorem 2.1, for $\frac{\varepsilon}{2} > 0,$ there exists
$\delta > 0$ such that for any two continuous maps $v, u :
S_{\varepsilon} \rightarrow  S_{\varepsilon}$ $\delta$--near and a
$\delta$--homotopy $j_{t} : S_{\varepsilon} \rightarrow
S_{\varepsilon}$ with $j_{0} = v_{|_{S_{\varepsilon}}}$ and $j_{1}
= u_{|_{S_{\varepsilon}}},$ there exists an
$\frac{\varepsilon}{2}$--homotopy $g_{t} : D \rightarrow
S_{\varepsilon}$ such that $g_{0}=v$ and $g_{1}=u$ and
$g_{t_{|_{S_{\varepsilon}}}}=j_{t}$ for all $t \in [0, 1].$ By
[3, p.40] for $\frac{\delta}{2} > 0,$ there exists $\lambda > 0$
such that if $\varphi, \psi : S_{\varepsilon} \rightarrow
S_{\varepsilon}$ are two continuous maps $\lambda$-near, then
$\varphi$ and  $\psi$ are $\frac{\delta}{2}$--homotopic.

Since $f_{\varepsilon}|_{S_{\varepsilon}}$ is a continuous map and
$S_{\varepsilon}$ is a finite polyhedron, then by Hopf's
construction [3, 8] there exists a continuous map $k: S_{\varepsilon}
\rightarrow S_{\varepsilon}$ which is $\lambda$-near to
$f_{\varepsilon}|_{S_{\varepsilon}}$ and has only a finite number
of fixed points. Hence, from [3, p.40], $k$ and
$f_{\varepsilon}|_{S_{\varepsilon}}$ are ${1\over
2}\delta$-homotopic. Let $(h_{t})_{( t\in [0,1])}$ be this
${1\over 2}\delta$-homotopy between $k$ and
$f_{\varepsilon}|_{S_{\varepsilon}}$ such that $h_{0}=f_{\varepsilon}|_{S_{\varepsilon}}$  and
$h_{1}=k.$ Now, we define a new homotopy
$(j_{t})_{ t\in [0,1]}$  by setting:
$$ j_{t}= \left\{ \begin{array}{ll}
                        h_{2t}  &\mbox{ if}\ 0\leq t<{1\over 2}\\
                        h_{2t-2} &\mbox{ if }\ {1\over 2}\leq t\leq 1.
                \end{array}
                \right. $$
So, we deduce that $(j_{t})_{t\in [0,1]}$ is a $\delta$-homotopy
such that $j_{0}=j_{1}=f_{\varepsilon}|_{S_{\varepsilon}}$ and
$j_{{1\over 2}}=k.$

Hence, by Theorem 2.1 there exists an  ${1\over
2}\varepsilon$-homotopy $g_{t}: D \rightarrow S_{\varepsilon}$
such that $g_{0}=g_{1}=f$ and $g_{t_{|_{S_{\varepsilon}}}}=j_{t}$
for all $t \in [0,1].$ Put, $g_{{1\over 2}}=g.$ Then, $g: D
\rightarrow S_{\varepsilon}$ is a continuous map and
$Fix(g)=Fix(k).$ Then, $Fix(g)$ is a finite set.  Since for every $x \in C$ we have
  $$||f(x) - g(x)|| \leq ||f(x) - f_{\varepsilon}(x)|| + ||f_{\varepsilon}(x) - g(x)||.$$
  So, we get $||f(x) - g(x)|| < \varepsilon.$ Thus, the map $g$ is
  $\varepsilon$--near to $f$ and has only a finite number of fixed
  points.

As a consequence of Theorem 3.1, we get the following result.

\begin{coro}
Let  $C$ be a nonempty compact convex  subset of an $F$-space
  $(E, || \quad ||),$ $f : C \rightarrow C$ be a continuous map   and $\varepsilon > 0$
be given. Then  there exists a continuous map $f_{\varepsilon} :
C \rightarrow C$ which is $\varepsilon$--near to $f$ and has only
a finite number of fixed points.
\end{coro}

\section{Fix--finite approximation  propriety for $n$--valued multifunctions}
In this section, By using  the simplicial approximation property (Lemma 2.7) and Theorem 3.1, we shall establish a fix--finite approximation  result for $n$--valued multifunctions in $F$--spaces. More precisely we shall show the following.

\begin{teo}
Let  $C$ be a nonempty compact convex  subset of an $F$--space
$(E, || \quad ||)$ and let $D$ be a path and simply connected compact subset
of $E$ containing $C.$ Then the pair $(D, C)$ satisfies
the F.F.A.P. for every continuous $n$--valued continuous multifunction $F : D \rightarrow 2^{C}.$
\end{teo}

In order to give the prove of Theorem 4.1, we shall need the following proposition.

\begin{prop}
Let  $C$ be a nonempty compact convex  subset of an $F$--space
$(E, || \quad ||)$ and let $D$ be a compact subset of $E$ containing $C.$
 Then the pair $(D, C)$ satisfies
the F.F.A.P. for any $n$--function $G : D \rightarrow 2^{C}.$
\end{prop}

\bigskip

\noindent Proof. Let $G : D \rightarrow 2^{C}$ be an $n$--function. Then, there exist $n$ continuous maps $g_{i} : D \rightarrow C$ such that $G(x)=\{g_{1}(x),...,  g_{n}(x)\}$ for all $x\in D$ and  $g_{i}(x) \not=g_{j}(x)$ for all $x \in D$ and $i,j=1,...,n$ with $i \not=j.$

For all $i,  j=1,...,n$ with $i \not=j,$ we define
$\gamma_{(i,  j)}(G)=\min\{ || g_{i}(x) - g_{j}(x) || : x \in D \}$. As every $g_{i}$ is continuous for all $i=1,..., n$ and $D$ is compact, then for
every $i,j=1,...,n$ with $i \not=j,$ we have $\gamma_{(i,  j)}(G)>0.$ Therefore,
 $$\gamma(G)=\min\{ \gamma_{(i,  j)}(G) : i,j=1,...,n, i \not=j \}>0.$$
 Let $\varepsilon>0$ be given. Then, we set
$\gamma=\min( {1\over 2} \gamma(G),  {1\over 2}\varepsilon ).$ Hence, $\gamma > 0.$
By Theorem 3.1, for each $i=1,..., n,$ there exists a continuous map
$h_{i} : D \rightarrow C$ which is $\gamma$-near to $g_{i}$ and has only a finite number of fixed points.  Let $H :  D\rightarrow C$ be the  multifunction defined by $H(x)=\{h_{1}(x),..., h_{n}(x)\},$ for all $x\in D.$ Then, the multifunction $H$ is an $n$--function. Indeed, if there exists
$x_{0}\in D$ and $i,j=1,...,n$ with $i \not=j,$ such that $h_{i}(x_{0})=h_{j}(x_{0})$,
then,
 $$|| g_{i}(x_{0}) - g_{j}(x_{0}) ||\leq || g_{i}(x_{0}) - h_{i}(x_{0}) || + || h_{j}(x_{0}) - g_{j}(x_{0}) || < 2\gamma.$$ So, we get $\gamma_{(i, j)}(G) < \gamma(G).$ This is a contradiction and $H$ is an $n$--function.  Now, since for all $i=1,...,n$ and for every $x \in D,$ we have,
 $|| g_{i}(x) - h_{i}(x) ||< {1\over 2}\varepsilon.$ So, we deduce that we have $d_{H}(G, H)< \varepsilon.$ Thus, $H$ is $\varepsilon$-near to $G.$ On the other hand, we know that $Fix(H)=\displaystyle\bigcup_{i=1}^{i=n}Fix(h_{i}).$
 As for all $i=1,...,n$ the set $Fix(h_{i})$  is  finite, hence we conclude that
 the multifunction $H$  has only a finite number of fixed points.

\bigskip

Now, we are ready to give the proof of Theorem 4.1.

\bigskip

\noindent Proof of Theorem 4.1. Let $\varepsilon > 0$ be given and $G : D \rightarrow 2^{C}$ be a $n$--valued continuous multifunction. Then from [13, p.76], we know that $\gamma (G) > 0.$ Set $\delta=\min({1\over 4}\varepsilon, {1\over 2}\gamma(G)).$ So, $\delta > 0.$ By Lemma 2.7, there exist a simplex $S_{\varepsilon}$ in $C$ and  a continuous map $h_{\varepsilon} : C \rightarrow
S_{\varepsilon}$ such that $||h_{\varepsilon}(x) - x || < \delta$ for every $x \in C.$
Now, we define a new continuous multifunction $H_{\varepsilon} : D \rightarrow C$ by setting: $$H_{\varepsilon}(x)=(h_{\varepsilon}\circ G)(x), \hbox{ for all }x \in D.$$
 Let $x \in D$ such that $G(x)=\left\{y_{1},..., y_{n} \right\},$ then $H_{\varepsilon}(x)=\left\{h_{\varepsilon}(y_{1}),..., h_{\varepsilon}(y_{n})\right\}.$ Assume that there is $i, j \in \{1,..., n\}$ such that $i\not=j$ and $h_{\varepsilon}(y_{i})=h_{\varepsilon}(y_{j}).$ So, we get
$$||y_{i} - y_{j}|| \leq ||y_{i} - h_{\varepsilon}(y_{i})|| + ||h_{\varepsilon}(y_{j}) - y_{j}||.$$
Thus, we have $||y_{i} - y_{j}|| < 2 \delta.$ Hence,  $||y_{i} - y_{j}|| < \gamma(G).$ Thus is a contradiction and
we deduce that the multifunction  $H_{\varepsilon}$ is $n$--valued. Now, since for all $i=1,..., n$ we have
$|| y_{i} - h_{\varepsilon}(y_{i}) ||<{1\over 4}\varepsilon,$ then we get $d_{H}(G, H_{\varepsilon})< {1\over 2} \varepsilon.$ Thus, $H_{\varepsilon}$ is ${1\over 2}\varepsilon$-near to $G.$ By Lemma 2.3, the multifunction $H_{\varepsilon} : D \rightarrow 2^{C}$ is an $n$-function. Using this and Proposition 4.2, we deduce that
there exists an $n$-function $K_{\varepsilon} : D \rightarrow 2^{C}$ which is ${1\over 2}\varepsilon$-near to $H_{\varepsilon}$ and has only a finite number of fixed points. As
$$d_{H}(G, K_{\varepsilon}) \leq d_{H}(G, H_{\varepsilon}) + d_{H}(H_{\varepsilon}, K_{\varepsilon}).$$
So, we get $d_{H}(G, K_{\varepsilon}) < \varepsilon.$ Therefore, we conclude that the multifunction
$K_{\varepsilon} : D \rightarrow 2^{C}$ is $\varepsilon$-near to $G$ and has only a finite number of fixed points.

As  consequences of Theorem 4.1, we get the following.

\begin{coro}
Let $C_{i}$  be a nonempty compact convex  subsets of an $F$--space $(E, || \quad ||)$
 for $i=1,..., n$  such that $\cap_{i=1}^{i=n} C_{i}\not=\emptyset.$ Let
 $k \in \{1, ..., n\}$ and let $G : \cup_{i=1}^{i=n} C_{i}\rightarrow C_{k}$  be a continuous $n$-valued multifunction and $\varepsilon > 0$ be given. Then, there exists a continuous $n$-valued multifunction $K_{\varepsilon} :
\cup_{i=1}^{i=n} C_{i} \rightarrow C_{k}$ which is $\varepsilon$--near to $G$ and has only
a finite number of fixed points.
\end{coro}

\begin{coro}
Let  $C$ be a nonempty compact convex  subset of an $F$--space
  $(E, || \quad ||),$ $G : C \rightarrow C$ be a continuous $n$--valued multifunction and $\varepsilon > 0$
be given. Then,

(i) $Fix(G)$ has at least $n$ fixed points ;

(ii)  there exists a continuous $n$--valued multifunction $K_{\varepsilon} :
C \rightarrow C$ which is $\varepsilon$--near to $G$ and has only
a finite number of fixed points.
\end{coro}

\bigskip

\noindent  Abdelkader Stouti \\
Center for Doctoral Studies: Sciences and Techniques   \\
Laboratory of Mathematics and Applications  \\
Faculty of Sciences and Techniques     \\
University Sultan Moulay Slimane    \\
P.O. Box 523 \\
Beni-Mellal 23000, Morocco  \\
E-mail address: stouti@yahoo.com or stout@fstbm.ac.ma

\end{document}